\documentclass[12pt]{article}

\title
{On the dynamic pull-in instability in  a mass-spring model of  
electrostatically actuated MEMS devices
}

\author{Gilberto Flores\\ 
Instituto de Investigaciones en Matem\'aticas 
Aplicadas y en Sistemas, \\ and FENOMEC\\ 
Universidad Nacional Aut\'onoma de 
M\'exico, \\Apdo. Postal 20-126, \\
01000 M\'exico, D.F., MEXICO  \\
(gfg@mym.iimas.unam.mx)}

\usepackage[dvips]{graphicx}
\usepackage{subcaption}

\def\complejo{\hbox{{\sf C}\kern-.45em\lower-.45ex\hbox{{\tiny ]}}\enskip}}
\def\real{\hbox{{\sf I}\kern-.1em{\sf R}}}
\def\natural{\hbox{{\sf I}\kern-.1em{\sf N}}}
\def\realc{\hbox{{\tiny I}\kern-.1em{\tiny R}}}
\def\entero{\hbox{{\sf Z}\kern-.40em{\sf Z}}}

\begin{document}

\baselineskip=.6cm

\maketitle

\vskip1cm

DEDICATED TO TIM MINZONI ON THE OCASSION OF HIS 65TH BIRTHDAY.

\vskip1cm

\newpage

\begin{abstract}
In this work we study the mass-spring system
\begin{equation}
\ddot x + \alpha \dot x + x =
- \frac{\lambda} {(1+x)^{2}},
\label{e:inertia}
\end{equation}
which is a simplified model for an electrostatically actuated MEMS device.
The static pull-in value is $\lambda^{*}=\frac{4}{27}$, which
corresponds to the largest value of $\lambda$ for which there
exists at least one stationary solution. For $\lambda > \lambda^{*}$ 
there are no stationary solutions and $x(t)$ achieves the value $-1$ in
finite time: {\it touchdown} occurs. We establish the existence
 of a dynamic pull-in value $\lambda_{d}^{*}(\alpha) \in (0, \lambda^{*})$, 
defined for $\alpha \in [0,\infty)$, which is a 
threshold in the sense that $x(t)$ approaches a stable stationary solution as 
$t \to \infty$
for $0 < \lambda < \lambda_{d}^{*}(\alpha)$, while touchdown occurs
for $\lambda > \lambda_{d}^{*}(\alpha)$. This dynamic pull-in value is a 
continuous, strictly increasing function of $\alpha$ and 
$\lim_{\alpha\to\infty} 
\lambda_{d}^{*}(\alpha)=  \lambda^{*}$.
\end{abstract}


\newpage

\section {\bf Introduction}

\vskip.3cm

The operation of many micro electromechanical systems (MEMS) relies upon the action
of electrostatic forces. Many such devices, including pumps, switches or valves,
can be modeled by electrostatically deflected elastic membranes. In a typical
situation, a MEMS device consists of an elastic membrane held at a constant voltage
and suspended above a rigid ground plate placed in series with a fixed voltage
source. The voltage difference causes a deflection of the membrane. 
For a more detailed description  we refer to the book by
Pelesko and Bernstein \cite{pelesko}.

Taking inertial and viscous forces into account, assuming that the membrane is
thin and using a linear approximation for the elastic energy, which is 
the analogue of a linear Hooke's law, the motion of the membrane is described
by a wave equation with damping and a singular forcing. Rescaling time yields,
 in the viscous 
dominated regime, the equation

\begin{equation}
\gamma^{2} u_{tt}+u_{t}-\Delta u = -\frac{\lambda}{(1+u)^{2}}\,\, {\rm in} \,\,
\Omega. 
\label{e:dampedwave}
\end{equation}

In the regime dominated by inertia the equation is:

\begin{equation} u_{tt}+\alpha u_{t}-\Delta u = -\frac{\lambda}{(1+u)^{2}} \,\, 
{\rm in}\,\,  \Omega. 
\label{e:inertialwave}
\end{equation}

where $\gamma$ is the ``{\it quality factor}'' and 
$\displaystyle \alpha=\frac{1}{\gamma}$.

An important nonlinear phenomenon in electrostatically deflected membranes is 
the so-called ``pull-in'' instability. 
For moderate values of the voltage the system is
in the {\it stable operation regime} in which the membrane approaches a stable 
steady state and remains separate from the ground plate; when the voltage  is 
increased beyond a critical value, the device is in the {\it touchdown regime}:
the membrane collapses onto the ground plate.
This phenomenon is known as ``touchdown'' or ``pull-in''. 

The critical value of the voltage required 
for touchdown to occur is termed the {\it pull-in voltage}. The determination
of the pull-in voltage is important for the design and manufacture of MEMS 
devices.
In most cases it is desirable to achieve the stable operation regime,
except for some devices such as microvalves, for which touchdown is a desirable
property.

Nathanson {\it et al.} \cite{nathanson} introduced the first model for an
electrostatically actuated  device, a millimiter-sized resonant gate
transistor was modeled by a mass-spring system.

In this model, the moving structure is a plate attached to a spring. 
The elastic properties of the
moving plate are described by the restoring force of the spring, which is 
assumed to be given by Hooke's law in the linear regime. 
The voltage applied to the moving plate results in an electrostatic force 
acting on the system by setting it in motion. 

The governing equation for the displacement of the moving mass 
is
\begin{equation}
\displaystyle m \ddot x + b \dot x + k x=-\frac{\lambda} {(1+x)^{2}},
\label{e:spring}
\end{equation}
in which the relevant parameter $\lambda$ is proportional to the square of the
applied voltage.

For systems in which damping dominates, we
introduce a dimensionless time with scaling factor
$k/b$, which yields

\begin{equation}
\displaystyle \gamma^{2} \ddot x + \dot x + x = 
- \frac{\lambda}{(1+x)^{2}},
\label{e:damping}
\end{equation}

where $\gamma^{2}=m k/ b^{2}$. In this formulation, it is easy to see 
that if the
inertia is not taken into account, that is, $\gamma=0$, then the static
and dynamic critical values coincide. Indeed, (\ref{e:damping})
 reduces to the first order equation 

\begin{equation}
\displaystyle \dot x + x + \frac{\lambda}{(1+x)^{2}}=0.
\label{e:forder}
\end{equation}

With $\displaystyle f(x,\lambda):= x+\frac{\lambda}{(1+x)^{2}}$ and
$\displaystyle g(x):=x(1+x)^{2}$, the stationary solutions of 
(\ref{e:forder}) are the solutions of $f(x,\lambda)=0$, which correspond 
to solutions of $g(x)=-\lambda$.
The cubic
$g$ has a local minimum at $x=-\frac{1}{3}$. The number of
stationary solutions in the region of interest $-1<x<0$ is determined by 
$\displaystyle \lambda^{*} := -g(-\frac{1}{3})=\frac{4}{27}$. 
There are two solutions $x_{1}(\lambda) < x_{2}(\lambda)$ for
$ 0<\lambda<\lambda^{*}$, one 
($\displaystyle x_{1}(\lambda)= x_{2}(\lambda)= -\frac{1}{3})$
at $\lambda=\lambda^{*}$ and none for $\lambda>\lambda^{*}$.

Moreover, for $ 0<\lambda<\lambda^{*}$ and 
for any initial condition in $(x_{1}(\lambda),\infty)$,  
the corresponding solution converges to $x_{2}(\lambda)$ as $t \to \infty$.

\vskip.3cm

At $\lambda=\lambda^{*}$ we have $x(t)\to -\frac{1}{3}$ as $t\to\infty$ provided
$x(0)>-\frac{1}{3}$, in particular if $x(0)=0$. For $\lambda>\lambda^{*}$ and
$x(0)>-1$, we have $x(t)\to -1$ and the value $-1$ is achieved in finite time.
This is known as {\it quenching} in the mathematical literature.

\vskip.3cm

The coincidence of the static and dynamic pull-in values has also been
established for the parabolic equation obtained by setting $\gamma=0$
in (\ref{e:dampedwave}). See for instance, Flores {\it et al.}
\cite{flores2}.

Based on numerical evidence, several authors have reported that the
the dynamic pull-in value is smaller than the static pull-in value,
both for the wave equation and for the mass-spring system. This means that
when inertia is taken into account, the moving structure may collapse onto the
susbstrate even if there is a stable stationary solution.
Chang and Levine \cite{chang} observed this behavior for the
conservative wave equation, which corresponds to
$\alpha=0$ in (\ref{e:inertialwave}), Kavallaris {\it et al.} 
\cite{kavallaris} in a nonlocal version of the conservative wave equation
and Flores \cite{flores} in the damped wave equation (\ref{e:dampedwave}).

\vskip.3cm

For the mass-spring system,
Zhang {\it et al.} \cite{zhang} described the dynamic pull-in as the collapse
of the moving structure towards the substrate, due to the combined action of
kinetic and potential energies. They also stated that, in general, {\bf 
dynamic pull-in requires a lower voltage to be triggered compared to the
static pull-in threshold}. 

The main result in the present work establishes the
existence of a dynamic pull-in value for the  mass-spring system.

The equation for systems dominated by inertia is obtained from (\ref{e:spring})
by introducing a dimensionless time scaled by the natural
frequency of the system $\omega = \sqrt{k/m}$, which
yields  (\ref{e:inertia}), with $\alpha = 1/{\gamma}$. See
Pelesko and Bernstein \cite{pelesko}. 
Our results are formulated for this regime. The existence of the dynamic 
pull-in value is obtained for all positive values of $\alpha$, so that the
viscous dominated regime also has this property.

We assume that the motion starts from rest: $x(0)=0=\dot x(0)$.

\vskip.3cm

We begin by writing 
(\ref{e:inertia}) as the first order system 

\begin{equation}
\dot x=y, \qquad \dot y = - [\alpha y + f(x,\lambda)]
\label{e:system}
\end{equation}

\vskip.3cm

The stationary solutions of (\ref{e:system}) are given by $y=0$,
$f(x,\lambda)=0$. The structure of the steady states is obtained from
the first order equation (\ref{e:forder}).
There are two solutions $(x_{1}(\lambda),0)$ and $(x_{2}(\lambda),0)$ for
$ 0<\lambda<\lambda^{*}$, one (with
$\displaystyle x_{1}(\lambda)= x_{2}(\lambda)= -\frac{1}{3})$
at $\lambda=\lambda^{*}$ and none for $\lambda>\lambda^{*}$.

\vskip.3cm

We now describe the main steps in the proof of the existence of a dynamic 
pull-in value. In Section 2 we determine explicitly the dynamic 
pull-in value for $\alpha=0$, namely $\lambda^{*}_{d}(0)=\frac{1}{8}$.
In Section 3 we prove that for each $\lambda\in (\frac{1}{8},\lambda^{*})$
there exists a unique value $\alpha^{*}(\lambda)$ of $\alpha$ such that
(\ref{e:system}) is in the touchdown regime for $\alpha<\alpha^{*}$ and it
is in the stable operation regime for $\alpha>\alpha^{*}$. This threshold
is a continuos, strictly increasing function of $\lambda$,  and
$\lim_{\lambda\to\lambda^{*}} \alpha^{*}(\lambda)=\infty$.
The dynamic pull-in value $\lambda^{*}_{d}(\alpha)$ is then the inverse
function of $\alpha^{*}(\lambda)$. In view of the asymptotic behavior
of $\alpha^{*}(\lambda)$ as $\lambda\to\lambda^{*} $, the
dynamic pull-in value is defined for all $\alpha>0$, 
$\lambda^{*}_{d}(\alpha)\in (\frac{1}{8},\lambda^{*})$
and $\lim_{\alpha\to\infty} \lambda^{*}_{d}(\alpha)=\lambda^{*}$. 
Therefore, our results remain valid in
the viscous dominated regime which corresponds to $\alpha$ large.
In the limiting case $\alpha=\infty$ we obtain the coincidence of the static
and dynamic pull-in values for the first order equation mentioned above. 
A key property in the analysis is the monotonicity of the stable manifold of
$(x_{1}(\lambda),0)$, which determines the domain of attraction of
$(x_{2}(\lambda),0)$.

We conclude this introduction by mentioning that
one of the findings in Rocha {\it et al.}
\cite{rocha} is the fact that for an overdamped device, the dynamics
in the touchdown regime has three distinguished regions charachterized
by different time scales:
in the first region the structure moves fast until it gets near the static
pull-in distance, then there is a metastable region of very slow motion
and finally a third region in which collapse takes place on a fast
time scale.
We shall see that these regions correspond to the approach to the unstable
stationary solution, which occurs on a fast
(order 1) time scale, followed by a slow motion close to the stable manifold
of the unstable steady state until the solution gets away from the stationary
point and enters the region of collapse where the dynamics occurs on a fast
time scale again.

\vskip.5cm

\section{The conservative case: $\alpha=0$.}

In this case it is possible to determine explicitly the 
dynamic pull-in value which separates the stable operation
regime from the touchdown regime. In the present situation, the
stable operation regime means that the solution is periodic.

\vskip.3cm

When $\alpha=0$,
(\ref{e:system}) becomes  a conservative system. 
The integral curves are determined explicitly as graphs of functions 
by means of

\begin{equation}
y=\pm \sqrt{2} \sqrt{E_{0} - F(x,\lambda)}
\label{e:integral}
\end{equation}

where $E_{0}$ is the total energy of an initial condition and
$\displaystyle F(x,\lambda)= \frac{x^{2}}{2} - \frac{\lambda}
{1+x}$ is a primitive of $f(x,\lambda)$.

For each $\lambda \in (0, \lambda^{*})$, $(x_{1}(\lambda),0)$ is a 
saddle, $(x_{2}(\lambda),0)$ is a center surrounded by periodic orbits
and a homoclinic orbit at $(x_{1}(\lambda),0)$.

It is clear from the picture in the phase plane that the solution starting
at $(0,0)$ is periodic if and only if this initial condition is enclosed
by the homoclinic. It is also clear that this happens if and only if
$F((x_{1}(\lambda),\lambda) > F(0,0)=-\lambda$. 
Using these observations, we prove 

\vskip.5cm

{\bf Proposition 1.}  
There exists $\lambda_{d}(0) \in (0, \lambda^{*})$ such that the solution
starting at $(0,0)$ is periodic if and only if $\lambda < \lambda_{d}(0)$,
and it approaches $(-1, -\infty)$ for $\lambda > \lambda_{d}(0)$.
In fact, $\lambda_{d}(0)=\frac{1}{8}$.

\vskip.3cm

{\sl Proof}.
Indeed, 
$\displaystyle \phi(\lambda):= F(x_{1}(\lambda),\lambda) - F(0,\lambda) =
\frac{x_{1}^{2}(\lambda)}{2} +\lambda \frac{x_{1}(\lambda)}{1+x_{1}(\lambda)}$
satisfies
$\displaystyle
\frac{d\phi}{d\lambda}= f(x_{1}(\lambda),\lambda) \frac{dx_{1}}{d\lambda}+ 
\frac{\partial F}{\partial \lambda}(x_{1}(\lambda),\lambda) - 
\frac{\partial F}{\partial \lambda}(0,\lambda) =
\frac{x_{1}(\lambda)}{1+x_{1}(\lambda)} <0.
$ 
Thus, $\phi$ is an strictly
decreasing function. Moreover,  
$\phi(\frac{4}{27})= -\frac{1}{54}$, 
and $ \lim_{\lambda \to 0^{+}} \phi(\lambda)= 1/2$ since
$\displaystyle \lambda \frac{x_{1}(\lambda)}{1+x_{1}(\lambda)} = 
- x_{1}^{2}(\lambda) [ 1+ x_{1}(\lambda)]$ and $x_{1}(\lambda) \to -1$ as
$ \lambda \to 0^{+}$. 
The continuity and monotonicity  of $\phi$ guarantee the existence
of a unique value $\displaystyle \lambda_{d}(0)\in (0,\frac{4}{27})$ 
such that $\phi(\lambda_{d}(0))=0$. This is the dynamic pull-in value.
The root is determined explicitly using the previous identity: 
$\displaystyle x_{1}(\lambda_{d}(0))= -\frac{1}{2}$ and 
$\displaystyle \lambda_{d}(0)=\frac{1}{8}$. In terms of the phase portrait,
this means that the homoclinic orbit at $(x_{1}(\lambda),0)$ crosses the
$x$-axis at a point $(\bar x(0),0)$ with $\bar x(0)>0$ for 
$\displaystyle \lambda \in (0,\frac{1}{8})$, while $\bar x(0)<0$ for 
$\displaystyle \lambda \in (\frac{1}{8}, \frac{4}{27})$. 
The required properties of the
solution starting at the origin follow from this.
The proof is finished.

\vskip.5cm

\section {The dissipative case: $\alpha >0$}.

We begin with the local stability analysis of the stationary solutions.

\vskip.3cm

It is clear  that 
$\displaystyle (-1)^{j} 
\frac{\partial f}{\partial x}(x_{j}(\lambda),\lambda)>0$
for $ 0<\lambda<\lambda^{*}$.
The jacobian matrix of the vector field at the stationary solution
$(x_{j}(\lambda), 0)$, which we denote by $A_{j}(\lambda,\alpha)$  
is given by 

$$ A_{j}(\lambda,\alpha)=\left(
\matrix{ 0 & 1 \cr
\displaystyle -\frac{\partial f}{\partial x}(x_{j}(\lambda),\lambda) & -\alpha}
\right).$$

Its characteristic polynomial is $\displaystyle
p(\mu)= \mu^{2} + \alpha \mu + 
\frac{\partial f}{\partial x}(x_{j}(\lambda),\lambda)$, with 
roots
\begin{equation}
\mu_{\pm} = - \frac{\alpha}{2} \pm 
\sqrt{\left(\frac{\alpha}{2}\right)^{2}- \frac{\partial f}{\partial x}(x_{j}(
\lambda),
\lambda)}
\label{e:evalues}
\end{equation}

It follows that for $\displaystyle \lambda \in (0,\frac{4}{27})$, 
$(x_{1}(\lambda),0)$ is a saddle, while
$(x_{2}(\lambda),0)$ is a stable node if 
$\displaystyle \alpha > 2 
\sqrt{\frac{\partial f}{\partial x}(x_{2}(\lambda),\lambda)}$,
and it is a stable focus for values of $\alpha$ such that the
reversed inequality holds. At $\displaystyle \lambda= \frac{4}{27}$
we have a degenerate stationary solution: 
$(x_{1}(\lambda),0) = (x_{2}(\lambda),0)$ with eigenvalues 
$\mu_{+} = 0$ and $\mu_{-}=- \alpha$.  

\vskip.3cm

The stable operation regime corresponds to the values of $\lambda$ for
which $x(t;\alpha, \lambda) \to x_{2}(\lambda)$ as $t \to \infty$.
In dynamical terms, this means that the initial condition $(0,0)$ 
belongs to the domain of attraction of $ (x_{2}(\lambda),0)$.

By means of a phase plane analysis, we establish the existence of
a dynamic pull-in value $\lambda_{d}^{*}<\lambda^{*}$ such that the stable 
operation regime is the interval $(0, \lambda_{d}^{*})$, while the 
touchdown regime corresponds to $(\lambda_{d}^{*}, \infty)$.

\vskip.3cm

Indeed, (\ref{e:system}) is a dissipative system with energy

\begin{equation}
E(x,y) = \frac{y^{2}}{2} + \frac{x^{2}}{2} - \frac{\lambda}{1+x}
\label{e:energy}
\end{equation}

such that along integral curves, 
$\displaystyle \frac{dE}{dt}= -\alpha y$.

It follows that the system does not have periodic or homoclinic orbits,
and  every solution which is bounded for $t\geq 0$ converges to a
stationary solution.

\vskip.3cm

We denote by $\gamma(t;\alpha,\lambda)$ the solution of (\ref{e:system})
with $\gamma(0;\alpha,\lambda)=(0,0)$, and
consider the relevant region of parameters: 
$\Omega=\{(\lambda,\alpha)\,: \, \lambda>0, \,\, \alpha\geq0\}$, which is
divided into $\Omega_{1}=\{(\lambda,\alpha)\in\Omega\, :\, 
\gamma(t;\alpha,\lambda)\to (x_{2}(\lambda),0)$ as $t\to\infty\}$,
$\Omega_{2}=\{(\lambda,\alpha)\in\Omega\, :\, 
\gamma(t;\alpha,\lambda)\to (x_{1}(\lambda),0)$ as $t\to\infty\}$,
and $\Omega_{3}=\{(\lambda,\alpha)\in\Omega\, :\, 
\gamma(t;\alpha,\lambda)\to (-1,-\infty)\,\, {\rm in\,  finite\,  time}\}$.
The set $\Omega_{1}$ corresponds to the stable operation regime,
$\Omega_{3}$ corresponds to the touchdown regime and $\Omega_{2}$
 corresponds to the critical behavior. We shall prove
that $\Omega=\cup_{j=1}^{3}\Omega_{j}$.

\vskip.3cm

Our first result concerning (\ref{e:system}) is the existence of a stable
operation regime for each $\alpha>0$. Indeed, we show  
that $\Omega_{1}$
contains a vertical strip in $\Omega$. We also give an explicit description
of part of the domain of attraction of the stable steady state. 
The energy $E$ defined in 
(\ref{e:energy}) and the euclidean distance $D(x,y)=x^{2}+y^{2}$
are useful tools in the analysis. 

\vskip.3cm

{\bf Proposition 2.}
For fixed $\alpha >0$ and $\lambda<\frac{1}{32}$, the set
$\displaystyle U=\{(x_{0},y_{0}): \quad 
D(x_{0},y_{0}) <\frac{1}{16}, \quad {\rm and} \quad
E(x_{0},y_{0})\leq -\lambda\}$ is positively invariant. Integral curves
of (\ref{e:system}) corresponding to initial conditions in $U$ satisfy
$x(t)\to x_{2}(\lambda)$ as $t \to \infty$.  

{\sl Proof.} For a given $(x_{0},y_{0}) \in U$, take $T>0$ such that 
$D(x(t),y(t))<\frac{1}{4}$ for $0\leq t\leq T$.
Since $x(t)\geq -\frac{1}{2}$ for $0\leq t\leq T$, it follows that 
$\displaystyle D(x(t),y(t))= 2 E(x(t),y(t)) +\frac{2\lambda}{1+x(t)}
\leq -2\lambda + 4 \lambda =2\lambda<\frac{1}{16}$. It follows that 
$(x(t),y(t)) \in U$ for $ 0\leq t\leq T$. Since $T$ depends on the
Lipschitz constant of the vector field on a fixed domain, we conclude
that  $(x(t),y(t)) \in U$ for all $t>0$. This establishes the positive
invariance of $U$. Moreover, $(0,0) \in U$. Therefore, the corresponding
integral curve converges to a stationary solution as $t \to \infty$.
In Proposition 1 we established that $x_{1}(\lambda)<-\frac{1}{2}$ for
$\lambda <\frac{1}{8}$. In particular, the same is true for the values
of $\lambda$ under consideration. Hence, $x(t)\to x_{2}(\lambda)$
as $t \to \infty$. The proof is finished.  

\vskip.3cm

The above result is similar to Theorem 2 of \cite{flores} in which
the existence of the stable operation regime is established for the
damped wave equation model. Since $(0,0)\in U$, it follows that
$\displaystyle(0,\frac{1}{32})\times (0,\infty)\subset \Omega_{1}$.

\vskip.3cm

The stable operation regime and the touchdown regime do persist under
small perturbations of the parameters and initial conditions. This 
implies that $\Omega_{1}$ and $\Omega_{3}$ are open subsets in $\Omega$.
This is the content of the following result.

\vskip.5cm

{\bf Proposition 3.} Denote by $\gamma(t;\alpha,\lambda)$ the integral curve
of (\ref{e:system}) with $\gamma(0;\alpha,\lambda)=(0,0)$. Assume that for
fixed $\alpha_{0}>0$ and $\lambda_{0}\in (0,\frac{4}{27})$ the corresponding 
integral curve satisfies either of the following two conditions:

a) 
$\gamma(t;\alpha_{0},\lambda_{0}) \to 
(x_{2}(\lambda_{0}),0)$ as $t \to \infty$ 

b) 
$\gamma(t;\alpha_{0},\lambda_{0}) \to (-1,-\infty)$ in finite time

Then the same is true for all nearby values of $\alpha$ and $\lambda$.

{\sl Proof}

a) Since $(x_{2}(\lambda_{0}),0)$ is a hyperbolic sink, Taylor's theorem
guarantees the existence of $\delta(\lambda_{0})>0$ such that
on the circle centered at this fixed point and radius $\delta$, 
the vector field defined by (\ref{e:system}) points inside 
the corresponding disk. By continuous dependence on parameters, the same
is true for all values of $\alpha$ and $\lambda$ sufficiently close to
$\alpha_{0}$ and $\lambda_{0}$ respectively. For such values of 
$\alpha$ and $\lambda$, the corresponding integral curve 
$\gamma(t;\alpha,\lambda)$ enters the invariant disk at some positive time
and it remains there for all later times. It follows that the integral
curve must approach $(x_{2}(\lambda),0)$, provided we restrict $\delta$
further if necessary, to make sure that the saddles 
$(x_{1}(\lambda),0)$ lie outside the invariant disk.

b) Since for any positive values of $\alpha$ and $\lambda$ we have
$\ddot x(0;\alpha,\lambda)=\dot y(0;\alpha,\lambda))=-\lambda$,
the integral curve starting at the origin enters the third quadrant
immediately and $\dot y(t)<0$ for sufficiently small positive values of $t$.
At some $t_{1}>0$ we must have $\dot y(t_{1})=0$ and the integral curve
enters the region $\dot y(t)>0$ immediately. By the hypothesis, the
integral curve cannot remain in this region for all $t>t_{1}$.
Therefore, there exists $t_{2}>t_{1}$ such that 
$\gamma(t_{2};\alpha_{0},\lambda_{0})$ 
is in the invariant region defined by $-1<x<x_{1}(\lambda_{0})$, 
$y<0$ and
$\dot y=-[\alpha_{0} y +f(x,\lambda_{0})]<0$. By the continuity of solutions
with respect to parameters, the same will be true for all nearby values
of $\alpha$ and $\lambda$. On each of the invariant regions corresponding
to such values of $\alpha$ and $\lambda$ we have $x(t)\to -1$ and
$y(t)\to -\infty$ in finite time. The proof is finished.  
 
\vskip.3cm

The argument in $(b)$ above allows us to show that 
$\gamma(t;\alpha,\lambda)$ either converges to a stationary solution or
else it approaches $(-1,-\infty)$.

\vskip.5cm

{\bf Proposition 4.} For any $(\alpha,\lambda)\in \Omega$, either
$\gamma(t;\alpha,\lambda)\to (x_{j}(\lambda),0)$ as $t \to\infty$ for
$j=1$ or $2$, or else $\gamma(t;\alpha,\lambda)\to (-1,\infty)$ in
finite time. In other words, $\Omega=\cup_{j=1}^{3} \Omega_{j}$.

{\sl Proof.} If  $\gamma(t;\alpha,\lambda)$ is bounded, then it is defined
for all $t>0$ and has a non-empty $\omega$-limit set, which  consists of 
stationary solutions. Since this set is finite, 
$\gamma(t;\alpha,\lambda)$ converges to a stationary solution as $t \to
\infty$.

To analyze the other case, assume that $\gamma(t;\alpha,\lambda)$ is unbounded.
The first observation is that the integral curve does not cross the line
$x=0$ for positive times, since $\displaystyle E(0,y)=\frac{y^{2}}{2} -
\lambda > -\lambda=E(0,0)$ for $y\neq0$. 
Therefore, $-1<x(t;\alpha,\lambda)<0$ for $t>0$ as long
as the solution is defined. It follows that the $y$ component  
is unbounded. Since it is bounded above by the maximum of $-f(x,\lambda)$,
it follows that the $y$ component must approach $-\infty$ through a 
sequence of times $t_{n}\to\infty$. This is possible only in the 
invariant region $-1<x<0$, $y<0$ and $\dot y<0$, where, as we have seen in
part (b) of the previous result, $x(t)\to -1$ and $y(t)\to-\infty$
in finite time.

\vskip.5cm

{\bf Remark}
A consequence of Proposition 4 is that $(\alpha,\lambda)\in\Omega_{3}$
for any $\alpha \geq 0$ and all $\lambda>\lambda^{*}$, since in this case
there are no stationary solutions.

\vskip.5cm

Our next result is part of the description of the phase portrait of
(\ref{e:system}), in which the invariant manifolds of the stationary solutions
play a fundamental role. We prove the existence of a heteroclinic orbit from 
$(x_{1}(\lambda),0)$ to $(x_{2}(\lambda),0)$ for 
$\displaystyle \alpha > 2 \sqrt{s(\lambda)}$, where 
$\displaystyle s(\lambda):= \sup \frac{f(x,\lambda)}{x-x_{2}(\lambda)}$.
Since $\displaystyle \frac{\partial f(x,\lambda)}{\partial x}=1 -
\frac{2\lambda}{(1+x)^{3}}$ is an increasing function of $x$, it follows
that 
$\displaystyle s(\lambda)=\frac{\partial 
f(x_{2}(\lambda),\lambda)}{\partial x}$.
Hence, there is a heteroclinic orbit for every value of $\alpha$
for which $(x_{2}(\lambda),0)$ is a stable node.
The nonlinearity in the equation is of the type of Fisher's 
equation, which explains the existence of the saddle-node connections
for the stated values of $\alpha$.

\vskip.5cm

{\bf Proposition 5.} 
For each 
$\displaystyle \alpha> 2\sqrt{s(\lambda)}$,
system (\ref{e:system}) has a heteroclinic connection from
$(x_{1}(\lambda),0)$ to $(x_{2}(\lambda),0)$.

{\sl Proof.}
We consider the triangular region defined by the line $x=x_{1}(\lambda)$,
with $y>0$, the segment $[x_{1}(\lambda), x_{2}(\lambda)]$, with $y=0$
and a line segment $y=m(x-x_{2}(\lambda))$ with 
$x_{1}(\lambda)\leq x\leq x_{2}(\lambda)$, and $m<0$. 
It is clear that on the horizontal and vertical sides of the triangle, 
the vector field defined by (\ref{e:system}) points inward the triangular 
region. We shall determine negative values of $m$ for which the vector
field also points inwards on the third side. Choosing ${\bf N}=(-m,1)$
as a normal vector for the slanted side of the triangle,
and denoting by ${\bf V}$ the vector field defined by (\ref{e:system}),
the condition on $m$ so that the vector field points inwards is
$<{\bf N}, {\bf V}> < 0$. Since
\begin{equation}
<{\bf N}, {\bf V}> =-[x-x_{2}(\lambda)]\{m^{2} +\alpha m + 
 \frac{f(x,\lambda))}{[x-x_{2}(\lambda)]}\}, 
\label{dotproduct}
\end{equation}

it follows that $<{\bf N}, {\bf V}> \leq -[x-x_{2}(\lambda)][m^{2} +
\alpha m + s(\lambda)]$. The quadratic polynomial in $m$ has roots
$m_{\pm} = -\frac{\alpha}{2} \pm \frac{1}{2} 
\sqrt{\alpha^{2} -4 s(\lambda)}$.
For $\alpha>2\sqrt{s(\lambda)}$, the root $m_{-}$ is negative. Therefore,
there are negative values of $m$ for which the quadratic takes
negative values. This completes the construction of an invariant region.
The branch of the unstable manifold of $(x_{1}(\lambda),0)$ that
points into the region $y>0$, $x>x_{1}(\lambda)$ enters this 
invariant region
and never leaves it. Therefore, it converges to $(x_{2}(\lambda),0)$
as $t \to \infty$. The proof is finished.

\vskip.5cm

The next result establishes the monotonocity of the integral curve
starting at $(0,0)$ as a function of $\lambda$, as well as a criterion
for touchdown. 

\vskip.5cm

{\bf Proposition 6.} Fix $\alpha>0$, let 
$\gamma(t;\alpha,\lambda)=(x(t;\alpha,\lambda),y(t;\alpha,\lambda)$, then
$y(t;\alpha,\lambda))$
is a decreasing function of $\lambda$ as long as $y(t;\alpha,\lambda)<0$. 
Moreover, $(\alpha,\lambda)\in \Omega_{3}$ 
if there exists $\lambda_{0}<\lambda$ such that
$(\alpha,\lambda_{0})\in \Omega_{2}\cup\Omega_{3}$.

{\sl Proof}.  Take $0<\lambda_{1}<\lambda_{2}$, and let 
$y_{j}(t;\alpha):=y(t;\alpha,\lambda_{j})$. Then $y_{j}(0;\alpha)=0$ and
$\dot y_{j}(0;\alpha)=-\lambda_{j}$. It follows that 
$y_{2}(t;\alpha)<y_{1}(t;\alpha)$ for $t>0$ and small.  

Note that  the second component of the vector field in (\ref{e:system}) 
is monotonic in $\lambda$ because
$\displaystyle \frac{\partial f(x,\lambda)}{\partial \lambda}=
\frac{1}{1+x^{2}}>0$. This implies that the inequality above is valid as long
as each $y_{j}$ is negative. This means that the integral curves do not cross
as long as they remain in the third quadrant.

For the second part of the
statement, the conclusion follows immediately from the monotonicity
if $(\alpha,\lambda_{0})\in \Omega_{3}$. In the other case,
$(\alpha,\lambda_{0})\in \Omega_{2}$,
the integral curve starting at $(0,0)$ approaches $(x_{1}(\lambda_{0}),0)$
as $t \to \infty$. Since $x_{1}(\lambda)$ is increasing in $\lambda$ and
$x_{2}(\lambda)$ is decreasing, the monotonicity of the 
integral curves guarantee that $\gamma(t;\alpha,\lambda)$ cannot approach either
of the critical points. Hence it has to approach $(-1,-\infty)$ and the 
integral curve has a finite time of existence since the $y$ component
is eventually decreasing.  The proof is finished.

\vskip.5cm

The stable manifold of the saddle $(x_{1}(\lambda),0)$ plays a crucial
role in the determination of the dynamic pull-in value.
The domain of attraction of $(x_{2}(\lambda),0)$ is determined by the
connected component of the stable manifold of $(x_{1}(\lambda),0)$ 
that approaches this saddle from the third quadrant. 
It is more convenient to analyze the behavior of the stable manifold
by fixing $\lambda$ and varying $\alpha$.

\vskip.3cm

We prove that for each $\lambda \in (0,\frac{4}{27})$,  
the connected component of the stable
manifold described above is a strictly monotonic function
of $\alpha$. The point of intersection  with the horizontal axis is
a monotonic and continuous function of $\alpha$. 
We shall consider $\lambda \in (\frac{1}{8}, \frac{4}{27})$. 
We shall prove that for small
values of $\alpha$, the stable manifold crosses the negative $x$-axis,
which corresponds to touchdown because the solution starting at $(0,0)$ 
cannot approach $(x_{2}(\lambda), 0)$. For large values of $\alpha$,
the stable manifold crosses the positive $x$-axis. In this case, the
solution $(x(t),y(t))$ is bounded for $t\geq 0$ and it converges to
$(x_{2}(\lambda), 0)$ as $t \to \infty$.

\vskip.3cm

It follows that there is a unique value
$\alpha^{*}(\lambda)$ of $\alpha$ 
such that the stable manifold crosses the $x$- axis at
$x=0$. 
We also prove that $\alpha^{*}(\lambda)$ is a continuous and strictly
increasing function of $\lambda$. The dynamic pull-in value 
$\lambda_{d}^{*}(\alpha)$ is the inverse function of $\alpha^{*}(\lambda)$.  
The dynamic pull-in value is defined for all(?) 
positive values of $\alpha$.

\vskip.3cm

For convenience, we change $t$
by $-t$,  $y$ by $-y$, and rewrite (\ref{e:system}) in terms of 
$u=x-x_{1}(\lambda)$ and $v=y$, obtaining 

\begin{equation}
\dot u=v, \qquad \dot v= \alpha v - f(u+x_{1}(\lambda), \lambda)
\label{e:systemb}
\end{equation}

This system has a saddle point at $(0,0)$,  with eigenvalues given by:

\begin{equation}
\mu_{\pm} =  \frac{\alpha}{2} \pm
\sqrt{\left(\frac{\alpha}{2}\right)^{2}- \frac{\partial f}{\partial x}
(x_{1}(\lambda),
\lambda)}
\label{e:evalues}
\end{equation}

The branch of the local unstable manifold that points into the first quadrant
is the graph of a continuous function $v=\Phi(u;\lambda,\alpha)$ and it can be 
continued as a graph as long as $v>0$. Moreover, $\Phi(0;\alpha,\lambda)=0$
and $\displaystyle \frac{d\Phi}{du}(0;\alpha,\lambda)=\mu_{+}$.
By the Chain Rule,

\begin{equation}
\frac{d\Phi}{du}=\alpha -\frac{f(u+x_{1}(\lambda),\lambda)}{\Phi(u;\alpha,
\lambda)}
\label{e:destabman}
\end{equation}

Our next result is the monotonicity of $\Phi$ with respect 
to $\alpha$.

\vskip.5cm

{\bf Proposition 7.}
For fixed $\lambda>0$, $\Phi(u;\alpha,\lambda)$ is an strictly increasing 
function of $\alpha$. 

{\sl Proof.}
Fix $\lambda >0$, take $0<\alpha_{1}<\alpha_{2}$ and denote $\Phi(u;\alpha_{j},
\lambda)$ by $\Phi_{j}(u)$ for $j=1, 2$. 
Since $\mu_{+}$ is an strictly increasing 
function of $\alpha$, it follows that $\Phi_{1}(u)<\Phi_{2}(u)$ for small
positive values of $u$. It is clear from (\ref{e:destabman}) that the
graph of $\Phi_{1}$ cannot intersect the graph of $\Phi_{2}$ as long as they
are defined. The proof is finished.

\vskip.5cm

For fixed $\lambda \in (\frac{1}{8}, \frac{4}{27})$ , let 
$$ I(\lambda):=\{ \alpha \geq 0:\,\, {\rm there \,\,exists}
\,\, \bar u(\alpha)>0 \,\, 
{\rm with}\,\,  \Phi(\bar u(\alpha); \alpha, \lambda)=0\}$$ 
and 
$\displaystyle J(\lambda) := \{\bar u(\alpha) : \alpha \in I(\lambda) \}$.

\vskip.3cm

A crucial step in the proof of the existence of the dynamic pull-in
value is the determination of the set $J(\lambda)$. To do this, it is
convenient to analyze the intersection of the unstable manifolds with
the vertical line ${\bf L}$ in the phase plane given by 
$u=x_{2}(\lambda)-x_{1}(\lambda)$. By the monotonicity of the unstable
manifolds and the transversality of ${\bf L}$ with respect to the vector
field in (\ref{e:systemb}), the set 
$$
K(\lambda):= \{(x_{2}(\lambda)-x_{1}(\lambda), 
\Phi(x_{2}(\lambda)-x_{1}(\lambda);\alpha,\lambda)\, : \,\alpha\geq 0\}
$$
defines an interval on the line ${\bf L}$, since the points of intersection
define a continuous function of $\alpha$. See Conley \cite{conley}.

\vskip.3cm

Let $v_{0}:=\Phi(x_{2}(\lambda)-x_{1}(\lambda);0, \lambda)$ denote the
height of the homoclinic orbit corresponding to $\alpha=0$ at
$u=x_{2}(\lambda)-x_{1}(\lambda)$.  Our next result determines the
set $K(\lambda)$.

\vskip.5cm

{\bf Lemma 1.} 
$\displaystyle K(\lambda)= \{x_{2}(\lambda)-x_{1}(\lambda)\} \times
 [v_{0}, \infty)$.

{\sl Proof}. Since $f(u+x_{1}(\lambda),\lambda)<0$, it follows from
(\ref{e:destabman}) that 
$\displaystyle \frac{d\Phi}{du}\geq \alpha$ for every
$\alpha >0$ and $u\in [0,x_{2}(\lambda)-x_{1}(\lambda)]$. On this interval
we have $\displaystyle \Phi(u;\alpha,\lambda)\geq \alpha u$. 
In particular 
$\displaystyle \Phi(x_{2}(\lambda)-x_{1}(\lambda);\alpha,\lambda)\geq \alpha
[x_{2}(\lambda)-x_{1}(\lambda)]$. The proof is finished.

\vskip.5cm

{\bf Remark.} The above estimate for $\Phi$ suggests that the critical
value of $\alpha$ tends to $\infty$ as $\lambda$ approaches $\lambda^{*}$.

\vskip.5cm

{\bf Lemma 2.} 
$I(\lambda)$ and $J(\lambda)$ are non-empty intervals. Moreover,
$I(\lambda)=[0, \bar\alpha(\lambda))$ with 
 $\bar \alpha(\lambda):=\,\,\sup\,\, I(\lambda)$, and
$\bar u(\lambda):= \sup J(\lambda)=\infty$, so that
$J(\lambda)=[\bar u(0),\infty)$.

{\sl Proof.}
In Section 2 we verified that $0\in I(\lambda)$. The monotonicity
and the continuity of
$\Phi$ with respect to $\alpha$ guarantee that $I(\lambda)$ and
$J(\lambda)$ are intervals.

We claim that $I(\lambda)=[0, \bar \alpha(\lambda))$.
This is clear if $\bar \alpha(\lambda)=\infty$. 
Now assume that the supremum is finite and that it belongs to the set 
$I(\lambda)$, then $\bar u(\lambda)<\infty$ and 
$\displaystyle \Phi(\bar u(\lambda);\bar \alpha(\lambda),\lambda)=0$. 
Take $\varepsilon>0$, then, by the continuity and monotonicity 
of $\Phi$ with respect to $\alpha$, the set
$$
\{(x_{2}(\lambda)-x_{1}(\lambda), \Phi(x_{2}(\lambda)-x_{1}(\lambda);\alpha,
\lambda):\alpha\in [0,\bar \alpha(\lambda) + \varepsilon]\}
$$ 
describes a closed interval   on the
line ${\bf L}$ in the phase plane.  
Moreover, the integral curve starting at $(\bar u(\bar\alpha(\lambda)),0)$
immediately enters the region $v<0$. By continuity with respect to
initial conditions, the same is true for values of $\alpha \in
(\bar\alpha(\lambda), \bar\alpha(\lambda)+\varepsilon)$ for 
$\varepsilon$ sufficiently small. 
This contradicts the definition of $\bar u(\lambda)$ as the
supremum of $J(\lambda)$ and  establishes the claim.

The next task is to verify that $\bar u(\lambda)=\infty$.
We distinguish two cases, according  to whether
$\bar \alpha(\lambda)$ is finite or infinity.

In the first case, $\bar \alpha(\lambda)<\infty$, 
assume that $\bar u(\lambda)<\infty$. From (\ref{e:destabman}) we see that
$\Phi$ has a finite derivative on every finite interval on which it is
defined. Using this fact and the continuous dependence on $\alpha$
we conclude that $\Phi(u;\bar \alpha(\lambda),\lambda)>0$ for all $u>0$.
By continuous dependence on initial conditions  we get unstable
manifolds for $\bar \alpha(\lambda)-\varepsilon < \alpha <\bar \alpha(\lambda)$
with $\varepsilon$ small enough, for which $\Phi>0$ for $u>\bar u(\lambda)$.
This is a contradiction since
$\Phi(\bar u(\alpha);\alpha, \lambda)=0$ for $0\leq \alpha <
\bar \alpha(\lambda)$. This contradiction proves that 
$\bar u(\lambda)=\infty$ if $\bar \alpha(\lambda)<\infty$.

In the case $\bar \alpha(\lambda)=\infty$ we have that $\bar u(\alpha)$
is defined for all $\alpha >0$ and $\Phi(\bar u(\alpha);\alpha, \lambda)=0$.
The point of intersection of the graph of $\Phi$ with the line ${\bf L}$
lies in the region where $\dot v>0$. In the case under consideration,
the graph of $\Phi$ must leave this region at a point with first component 
$u$ satisfying $u>h(\alpha[x_{2}(\lambda)-x_{1}(\lambda)], \lambda)$,
where $h$ is the inverse function of $f$ on the interval 
$(x_{2}(\lambda)-x_{1}(\lambda),\infty)$. Since $f$ is increasing and
tends to $\infty$ as $u$ approaches $\infty$, we get 
$\displaystyle \bar u(\alpha)>h(\alpha[x_{2}(\lambda)-x_{1}(\lambda)],\lambda)
\to \infty$ as $\alpha \to \infty$. The proof is finished.

\vskip.5cm

For the values of $\lambda$ under consideration, the next result establishes
the existence of a critical value $\alpha^{*}(\lambda)$ such that the
touchdown regime corresponds to $(0, \alpha^{*}(\lambda))$, while the
stable operation corresponds to $(\alpha^{*}(\lambda),\infty)$.
The critical value occurs when $\bar u(\alpha)=  - x_{1}(\lambda)$.

\vskip.3cm

{\bf Theorem.} 
For each $\lambda \in (\frac{1}{8},\frac{4}{27})$ 
there exists $\alpha^{*}(\lambda)>0$
such that the touchdown regime corresponds to $(0, \alpha^{*}(\lambda))$, 
while the stable operation corresponds to $(\alpha^{*}(\lambda),\infty)$.
Moreover, 
$\alpha^{*}(\lambda)$ is an strictly increasing, continuous function of
$\lambda$, and $\alpha^{*}(\lambda) \to \infty$ as $\lambda \to 
\frac{4}{27}$.

{\sl Proof.}  The critical value $\alpha^{*}(\lambda)$ is the value of
$\alpha$ for which $\bar u(\alpha)=  - x_{1}(\lambda)$. It is well defined
by Lemma 2. The monotonicity is a consequence of Proposition 7.
To verify the stated properties, we return to the original equation
(\ref{e:system}). In this setting, the critical value satisfies
$\bar x(\alpha^{*}(\lambda))=0$.

The content of Lemma 2 is that for $\alpha < \bar\alpha(\lambda)$,
the branch of the stable manifold of the saddle 
$(x_{1}(\lambda),0)$ that enters from the third
quadrant intersects the horizontal axis at $(\bar x(\alpha),0)$
and the points of intersection comprise the interval $[\bar u(0),\infty)$,
or equivalently, the interval $[\bar x(0)-x_{1}(\lambda),\infty)$
where $\bar u(0)<-x_{1}(\lambda)$, or equivalently
$\bar x(0)<0$. The left end-point of the interval is thus determined
by the homoclinic orbit in the conservative case $\alpha=0$. 
By the monotonicity and continuity of the points of
intersection, there exists a unique value $\alpha^{*}$ of $\alpha$ such
that the point of intersection satisfies  $\bar x(\alpha^{*})=0$.
For $\alpha <\alpha^{*}$, the point of intersection satisfies  
$\bar x(\alpha)<0$. In this case, it is clear that $(0,0)$ is
not in the domain of attraction of $(x_{2}(\lambda),0)$. It follows that
the integral curve starting at $(0,0)$, which enters the third quadrant 
immediately, in fact, it enters the region where $y<0$
and $\dot y<0$. Since  integral curves in this region cannot approach the point
$(x_{2}(\lambda),0)$, it follows that there exists $T_{1}>0$ such that
$\dot y(T_{1})=0$ and the integral curve enters immediately the region where
$\dot y>0$. But, it cannot remain there for all $t\geq T_{1}$
since it cannot cross the stable manifold. It follows that there exists
$T_{2}>T_{1}$ such that $\dot y(T_{2})=0$. Now the integral curves enters
the region where $-1<x<x_{1}(\lambda)$, $y<0$ and $\dot y<0$.
As it was established in Proposition 4, solutions in this positively 
invariant region satisfy
$x(t)=-1$ is achieved in finite time.

In the case $\alpha >\alpha^{*}$, the point of intersection satisfies 
$\bar x(\alpha)>0$.
The stable manifold provides a lower bound on the second component of the
integral curve starting at $(0,0)$. The graph $v=-\frac{1}{\alpha}
f(x,\lambda)$ on the interval $[x_{1}(\lambda),x_{2}(\lambda)]$ provides
the upper bound. Hence, the integral curve converges to 
$(x_{2}(\lambda),0)$ as $t \to \infty$.

For the continuity, take $\lambda_{0}$ in the interval under consideration,
and $0< \varepsilon<\alpha^{*}(\lambda_{0})$, then 
$\gamma(t;\alpha^{*}(\lambda_{0})-\varepsilon,\lambda_{0})
\to (-1, -\infty)$ 
and $\gamma(t;\alpha^{*}(\lambda_{0})+\varepsilon,\lambda_{0})
\to (x_{2}(\lambda),0)$. By Proposition 3, there exists $\delta >0$
such that the above properties are mantained if $|\lambda-\lambda_{0}|<
\delta$. Now we use the continuity and monotonicity of the unstable
manifolds to conclude that for such values of $\lambda$ we have
$\alpha^{*}(\lambda_{0})-\varepsilon < \alpha^{*}(\lambda)<
\alpha^{*}(\lambda_{0})+\varepsilon$, which establishes the continuity
of the critical value $\alpha^{*}(\lambda)$ 

The last step is the asymptotic behavior of $\alpha^{*}(\lambda)$
as $\lambda\to\lambda^{*}$. The proof is by 
contradiction. Assume that for $\displaystyle\lambda=\lambda^{*}=
\frac{4}{27}$, there exists a positive real number $\alpha_{0}$ such that
if $\gamma(t;\alpha)$ is the solution of (\ref{e:systemb}) for 
$\displaystyle\lambda=\lambda^{*}$ with $\gamma(0;\alpha)=(0,0)$, 
except that we keep the variables $(x,y)$, then
$\displaystyle\gamma(t;\alpha_{0})\to (-\frac{1}{3},0)$ as $t\to\infty$.
In this case, $(\alpha,\lambda)$ is in the stable operation regime for
all $\alpha>\alpha_{0}$ and $0<\lambda<\lambda_{*}$. Since 
$(\alpha,\lambda)$ is in the touchdown regime for all $\alpha>0$ and 
$\lambda>\lambda_{*}$, it follows that 
$\displaystyle\gamma(t;\alpha)\to (-\frac{1}{3},0)$ as $t\to\infty$
for all $\alpha>\alpha_{0}$. Now we have a one-parameter family of unstable
manifolds of $(-\frac{1}{3},0)$ with a branch that points into the second
quadrant and such that their first crossing with the line $y=0$ occurs at
$x=0$. Each of these branches of the unstable manifolds is the graph
of a function $y=\Phi(x;\alpha)$ defined for $-\frac{1}{3}\leq x\leq 0$.
Moreover, $\displaystyle \frac{d\Phi}{dx}(x;\alpha)=\alpha-
\frac{f(x;\lambda^{*})}{\Phi(x;\alpha)}$. In particular
$\displaystyle \frac{d\Phi}{dx}(-\frac{1}{3},\alpha)=\mu_{+}=\alpha$.
An argument similar to the one used in Proposition 7 for $0<\lambda<
\lambda_{*}$ shows that $\Phi$ is an increasing function of $\alpha$
for $-\frac{1}{3}\leq x\leq 0$. Now let $(x_{0},y_{0})$ be the point of
intersection of the branch of the unstable manifold with the graph of
$y=f(x,\lambda_{*})$, then $y_{0}=\Phi(x_{0},\alpha_{0})$. Since $\Phi$
is an increasing function of $x\in (-\frac{1}{3},x_{0})$ for
$\alpha\geq\alpha_{0}$,  we take $\delta>0$ such that 
$\displaystyle \Phi(x_{0}-\delta,\alpha_{0})=\frac{y_{0}}{2}$, 
and  consider values
of $\displaystyle\alpha>\frac{4\lambda_{*}}{y_{0}}$. 
Then, for $x\in[x_{0}-\delta,x_{0}]$
we have
$$\frac{d\Phi}{dx}(x;\alpha)=\alpha-\frac{f(x;\lambda^{*})}{\Phi(x;\alpha)}
\geq \alpha-\frac{2\lambda^{*}}{y_{0}} > \frac{\alpha}{2}.$$

It follows from the Mean Value Theorem that
$$\Phi(x;\alpha)= \Phi(x_{0}-\delta;\alpha) +\Phi(x;\alpha)-
\Phi(x_{0}-\delta;\alpha)\geq \frac{y_{0}}{2}+
\frac{\alpha}{2}(x-(x_{0}-\delta))
$$
and at $x=x_{0}$ we get $\displaystyle\Phi(x_{0};\alpha) > \frac{y_{0}}{2}+
\frac{\alpha}{2}\delta$. Now we take $\alpha$ such that 
$\displaystyle\frac{y_{0}}{2}+\frac{\alpha}{2}\delta>\lambda^{*}=
f(0,\lambda^{*})$,
then for such values of $\alpha$ the following inequality holds:
$\Phi(0;\alpha)>\lambda^{*}>0$. This contradiction shows that the assumption
is not valid. The proof is finished. 

\vskip.5cm

{\bf Acknowledgements}. I am deeply grateful to Arturo Olvera for several 
stimulating conversations on this problem and for several useful suggestions.

\end{document}